\documentclass[reqno,final,11pt]{amsart}
\usepackage{natbib}  
\usepackage{fancyhdr} 
\usepackage{color} 
\usepackage{hyperref} 
\usepackage{graphicx} 
\usepackage{geometry}
\usepackage{comment}
\usepackage[utf8]{inputenc}
\usepackage[T1]{fontenc}
\usepackage{verbatim}

\usepackage{tikz}
\usetikzlibrary{calc,shapes,backgrounds,arrows}
\usepackage{amssymb}
\include{styleALEA}

\definecolor{aleacolor}{rgb}{0.16,0.59,0.78}

\hypersetup{
breaklinks,
colorlinks=true,
linkcolor=aleacolor,
urlcolor=aleacolor,
citecolor=aleacolor}

\renewcommand{\cite}{\citet}

\theoremstyle{plain}
\newtheorem{theorem}{Theorem}[section]

\theoremstyle{definition}

\theoremstyle{remark}
\newtheorem{remark}[theorem]{Remark}

\makeatletter \@addtoreset{equation}{section} \makeatother

\begin{document}

\title[Berry--Esseen Theorem for Sample Quantiles]{Berry--Esseen Theorem for Sample Quantiles with Locally Dependent Data}

\author{Partha S.~Dey}
\address{Department of Mathematics, University of Illinois at Urbana-Champaign, IL, USA}
\email{psdey@illinois.edu} 

\author{Grigory Terlov}
\address{Department of Statistics and Operational Research, University of North Carolina at Chapel Hill, NC, USA}
\email{gterlov@unc.edu} 


\begin{abstract}
  We derive a Gaussian Central Limit Theorem for the sample quantiles based on locally dependent random variables with explicit convergence rate. Our approach is based on converting the problem to a sum of indicator random variables, applying Stein's method for local dependence, and bounding the distance between two normal distributions. We also generalize this approach to the joint convergence of sample quantiles with an explicit convergence rate.
\end{abstract}

\maketitle

\section{Introduction}

The Central Limit Theorem (CLT) is one of the fundamental theorems in probability theory and statistics. In classical form, it states that the sum of i.i.d.~finite variance random variables, appropriately centered and scaled, converges in distribution to the standard normal distribution. Since then, CLT has been strengthened and extended to various settings. The study of the asymptotic distribution of the sample median was first developed for the case of a continuous random variable by Sheppard in 1890 (according to~\cite{Hald98}). Asymptotic properties of sample quantiles of independent random variables have been extensively studied for both continuous distributions (see \eg~\cite{DavidNagaraja03} and references therein) and discrete ones (see  \eg~\cite{Ma11}). The CLT for the sample median of i.i.d.~random variables is derived either by converting the problem to a sum of indicators~\cite{Pollard}, by the Delta method~\cite{vdvbook}, or through Bahadur representation~\cite{Bahadur66} (known as Bahadur--Kiefer theorem due to its generalization~\cite{Kiefer67}). The last one is particularly useful as it gives an almost sure representation of the empirical quantile as a sum of i.i.d.~random variables and a further error term of order $n^{-3/4}\log n$. This, in turn, reduces the question of CLT and the rate of convergence to the classical setting. This method also has been generalized to tackle some of the dependent cases, which include $m$-dependent sequences~\cite{Sen68}, stationary sequences~\cite{Dembinska14,Hesse90}, $\phi$-mixing process~\cite{Sen72}, and auto-regressive processes~\cite{Dutta71} among others. While exact representations given by Bahadur--Kiefer type theorems yield Gaussian convergences in all these cases, the rates are not immediate, if achievable at all, as in the independent case.

In this work, we build on the elementary technique of converting sample quantiles to a sum of indicators as in~\cite{Pollard} and combine it with Stein's method to derive a Gaussian convergence for sample qualities with local dependencies. The class of dependencies that we consider is a direct generalization of the $m$-dependent sequences, where we require to know only some properties of dependency neighborhoods (see Assumption~\ref{as:dn}). Hence, our result can be applied to models not only for linear dependency structures but also for more graphical ones. Another benefit of our approach is that it gives an explicit rate of convergence and directly generalizes to a multivariate version. We also discuss the optimality of the derived rate of convergence in Section~\ref{sec:optimal}.

\subsection{Stein's method} Stein's method bounds the distance between the random variable of interest $W$ and the standard normal variable $Z$ in the following way
\begin{align}\label{eq:Steineq}
	d(W,Z)\le\sup_{f\in\cD}\abs{\E(f'(W)-Wf(W))}
\end{align}
for an adequately chosen class of functions $\cD$ depending on the metric $d(\cdot,\cdot)$. When $W=\sum_{i=1}^nX_i$, it turns out that it is often easier to work with the right-hand side of~\eqref{eq:Steineq}, even if $ X_i$ have dependencies between each other. Depending on the structure of such dependencies, one would bound $\E(f'(W)-Wf(W))$ in different ways. Hence Stein's method can be used with variety of approaches such as exchangeable pairs~\cite{Stein72}, dependency neighborhoods or local dependencies~\cite{Baldi89,ChenShao04,Fang16,FangKoike21, Rinott94}, size-bias~\cite{GoldsteinRinott96} and zero-bias couplings~\cite{Goldstein22,GoldsteinReinert96}, Stein coupling~\cite{ChenRollin10}, and through Malliavin calculus~\cite{NourdinPeccati10} among others. In Section~\ref{ssec:LD}, we briefly discuss the basics of Stein's method for local dependence (also known as the dependency neighborhoods approach) and refer to~\cite{ChenGoldstein11, DiaconisHolmes98, Ross11} for further reading on the topic.

One particular advantage of Stein's method is that $\E(f'(W)-Wf(W))$ naturally preserves the additive structure of random variables. Hence, when $W=\sum_{i=1}^nX_i$, Stein's method allows working with local interactions among $X_i$'s to derive global convergence. On the other hand, when one works with non-additive functions, such as sample quantiles, applying Stein's method requires additional effort.

\subsection{Notations}\label{ssec:notations}

Throughout this paper, we will use the following notations
\begin{itemize}

	\item $Z$ always denotes a standard normal random variable, $\Phi,\phi$ denote the distribution and density function of $Z$, respectively.
	\item $\mvZ=(Z_1,Z_2,\ldots,Z_\ell)$ denotes an $\ell$-dimensional standard Gaussian vector.
	\item $\cL(W)$ - the law of random variable $W$,
	\item $\overline{X}:=X-\E X$ denotes a centered version of a random variable $X$.
	\item $W_i$ represents the $i^{\textrm{th}}$ coordinate of a vector $\mvW$,
	\item $\dks(W,Z):=\sup_{x\in\bR}\abs{\pr(W\leq x)-\pr(Z\leq x)}$ -- Kolmogorov-Smirnov distance,
	\item $\dwas(W,Z):=\dwas(\cL(W),\cL(Z)):=\sup_{h: 1\textrm{-Lip.}}\abs{\E h(W) -\E h(Z)}$ - Wasserstein distance,
	\item $\dtv(W,Z):=\dtv(\cL(W),\cL(Z)):=\sup_{\text{Borel set } A}\abs{\pr(W\in A) -\pr (Z\in A)}$ - Total variation distance.

	\item  $f\lesssim g$ if $f=O(g)$, $f\ll g$ if $f=o(g)$ and $f\approx g$ if $c'g\le f\le cg$ for some universal constants $c,c'>0$,
	\item $x\vee y:=\max\{x,y\}$ and $x\wedge y:=\min\{x,y\}$.
\end{itemize}
Recall that for a matrix $V=(\!(v_{ij})\!)_{1\le i\le n,1\le j\le m}$ the Frobenius norm of $V$ is defined as
\[
	\norm{V}_{\F}:=\sqrt{\sum_{i,j}\abs{v_{ij}}^2}
\]
and the operator norm of $V$ is defined as
\[
	\norm{V}_{\op}:=\sup\{\norm{V\vu}_2/\norm{\vu}_2,\vu\neq0\}.
\]
\subsection{Setup}
Let $X_1, X_2, \ldots, X_n$ be real-valued random variables (not necessarily independent) with CDF's $\{F_i(x)\}_{i\in[n]}$, respectively. Given $\ga\in(0,1)$, we define the $\alpha^{\textrm{th}}$ sample quantile as
\begin{align}\label{def:Q}
	Q_{n}^{(\ga)}:=\sup\left\{x\in\dR\mid\abs{\{i\,|\,X_i\le x\}}\le\lfloor n\ga\rfloor\right\}.
\end{align}
In particular, the sample median is defined as
\begin{align}\label{def:M}
	M_n=Q_{n}^{(\half)}:=\sup\left\{x\in\dR\mid\abs{\{i\,|\,X_i\le x\}}\le\lfloor n/2\rfloor\right\}.
\end{align}

Fix an integer $\ell\ge 1$ and two increasing sequences of real numbers $0<\ga_1<\ga_2<\cdots<\ga_{\ell}<1$ and $m_{1}<m_{2}<\cdots<m_{\ell}$.
We assume the following.
\begin{ass}\label{as1}

	\begin{enumerate}
		\item For all $k\in[\ell], i\in[n]$, the equations $F_i(x)=\alpha_k$ has a unique solution at $x=m_{\ga_k}=:m_k$,\label{as: unique med q}
		\item For all $k\in[\ell], i\in[n]$, the CDF $F_i$ is continuously twice differentiable at $m_k$,\label{as: cont q}
		\item There exists $\eps, A>0$ such that for all $\abs{y}\le\eps$ we have\label{as:unifbdd}
		\begin{equation}\label{eq:unifbdd}
		\abs{F_i''(m_k+y)},\abs{F_i'(m_k+y)}\le A\text{ for all } k\in[\ell],\,i\in[n].
		\end{equation}
	\end{enumerate}
\end{ass}
We define
\begin{align}\label{def:c}
	\begin{split}
		\mvmu&:=(m_1,m_2,\cdots,m_\ell),\\
		\theta_{n}(x)&:=\frac{1}{n}\sum_{j=1}^nF'_j(x)\quad \textnormal{for } x\in\dR,\\
		\text{ and }\quad\Theta_{n}&:=\diag \left(\theta_{n}(m_1),\theta_{n}(m_2),\ldots,\theta_{n}(m_\ell)\right).
	\end{split}
\end{align}
When $\ga=1/2$, for simplicity, we may assume that $m_{\ga}=0$, and thus
\begin{align}\label{def:cmed}
	\theta_n:=\theta_n(0):=\frac{1}{n}\sum_{j=1}^nF'_j(0).
\end{align}

In this article, we consider the case when $(X_i)_{i\in[n]}$ is a locally dependent sequence of random variables. In particular, we will quantify the dependence via dependency neighborhoods. Dependency neighborhoods are present in various forms in the literature, see~\cite{ BarbourKaronski89, ChenShao04, Fang16,   paulin2012concentration,RinottRotar96} among others. For instance, one of the most common notions of generating dependency neighborhoods is via dependency graphs.
\begin{defn}[Dependency graph]\label{def:depgraph}
    A graph $G=([n],E)$ is called a dependency graph for random variables $(X_i)_{i\in[n]}$ if $(X_i)_{i\in A}$ and $(X_j)_{j\in B}$ are mutually independent whenever the sets of vertices $A\subseteq [n]$ and $B\subseteq [n]$ share no edges. 
\end{defn}

In this paper, we consider a more general notion of dependency neighborhoods given in the following assumption.

\begin{ass}\label{as:dn}
	We assume that the sequence of random variables $(X_1, X_2,\ldots, X_n)$ has the following structure of dependency neighborhoods, \ie
	\begin{enumerate}
		\item $\forall i\in[n]$, $\exists N_i\subseteq[n]$ such that $X_i$ is independent of $(X_l)_{l\in N_i^c}$.
		\item $\forall i\in[n]$ and $j\in N_i$, $\exists N_{ij}\subseteq[n]$ such that $N_{i}\subset N_{ij}$ and $(X_i,X_j)$ are independent of $(X_k)_{k\in N_{ij}^c}$,
		\item $\forall i\in[n]$, $j\in N_i$, and $k\in N_{ij}$, $\exists N_{ijk}\subseteq[n]$ such that $N_{ij}\subset N_{ijk}$ and $(X_i,X_j,X_k)$ are independent of $(X_s)_{s\in N_{ijk}^c}$.
	\end{enumerate}
	
\end{ass}

We refer to the set $N_i=\{j\mid j\sim i\}\cup\{i\}$ as the dependency neighborhood of $X_i$ and let
	\begin{equation}\label{eq:D}
	D_1:=\max_i\abs{N_i},\quad         
        D_2:=\max_{i,j}\abs{N_{ij}},\quad\text{and}\quad D_3:=\max_{i,j,k}\abs{N_{ijk}}.
	\end{equation}
These quantities are \textit{the main parameters of a dependency neighborhood and may depend on $n$}. Notice that the dependency graph condition is stronger than Assumption~\ref{as:dn} as it provides a natural choice for the dependency neighborhoods such that $D_2\le D_1^2$ and $D_3\le D_1^3$.

Given a vector $\mvx\in\dR^{\ell}$, we consider the centered random vector
\begin{align}\label{eq:vec Y}
	\mvY_{i,\mvx}:=\left(\1_{X_i\le m_k+ n^{-\half}\cdot x_k} - F_i(m_k+ n^{-\half}\cdot x_k)\right)_{k=1}^\ell.
\end{align}

Notice that the random vectors $(\mvY_{i,\mvx})_{i=1}^n$ inherit the structure of dependency neighborhoods from $( X_i)_{i=1}^n$.
Let
\[
	\gS_{n,\mvx} =\frac1n\var\left(\sum_{i=1}^n\mvY_{i,\mvx}\right)
\]
be the variance-covariance matrix of the random vector $\frac{1}{\sqrt{n}}\sum_{i=1}^n\mvY_{i,\mvx}$. In the univariate case, we denote it as $\gs_{n,x}^2.$
Finally, to simplify notation, we denote 
\begin{equation}\label{def:sigmas}
    \gs_n:=\gs_{n,0}\quad \text{and} \quad \gS_n:=\gS_{n,\mvzero}.
\end{equation} 
We will assume that the matrix $\gS$ is invertible.

\begin{remark}\label{rem:gs indep}
	In particular, if $\ga=\half$ and $m_\ga=0$ and in addition to the Assumption~\ref{as1} we know that $
		\pr(X_i\vee X_j\le 0)=\frac14$ for all $i\neq j,$
	then $\1_{X_i\le0}$ and $\1_{X_j\le0}$ are uncorrelated for $i\neq j$ and hence $\gs\equiv\half$. Similarly, when the random variables are independent, we have
	\[
	\gS_{n}=(\!(\,\ga_{i\wedge j} - \ga_i\ga_j\,)\!)_{i,j=1}^{\ell},
	\]
    which does not depend on $n$.
\end{remark}


\subsection{Main results}
In this section, we present our main results covering the univariate case (see Theorem~\ref{theorem: univ DG med}) and its extension to the multivariate case (see Theorem~\ref{theorem: multi quantile gen}). To present the univariate result, we assume that $\ga=\half$ and $m_\ga=0$; however, our argument can be easily extended to any other quantile that satisfies Assumption~\ref{as1}. To keep the exposition simple, we make a natural assumption that some of the parameters are uniformly bounded~\eqref{eq:unifbdd_thm} and~\eqref{eq:multbdd}, which could be relaxed by changing the Theorems~\ref{theorem: univ DG med} and~\ref{theorem: multi quantile gen} to hold for $n$ large enough, see Remark~\ref{rem:as}.

\begin{theorem}\label{theorem: univ DG med}
	Let $(X_i)_{i\in[n]}$ be a sequence of random variables satisfying Assumption~\ref{as1} with $\ga=\half$, $m_\ga=0$ and Assumption~\ref{as:dn}. Let $M_n$ be as in~\eqref{def:M}, $A$ be as in \eqref{eq:unifbdd}, $\theta_n$ be as in~\eqref{def:cmed},  $D_i$ as in \eqref{eq:D}, $\gs_n$ as in \eqref{def:sigmas}, and $Z\sim\N\left(0,1\right)$. Assume that there are constants $0<c<C<\infty$ such that for all $n\in\bN$
    \begin{equation}\label{eq:unifbdd_thm}
        c<\gs_n^2,\theta_n<C\quad\text{ and }\quad \max\{D_1,D_2,D_3\}\le C\cdot n^{\frac{1}{4}}.
    \end{equation}
    Then 

    \[
    \dks\left(\frac{\theta_n\sqrt{n}}{\gs_n }M_n\,,\,Z\right)\lesssim D_1\sqrt{\log n + D_1} \cdot\sqrt{\frac{\log n}{n}}+\frac{\max\{D_1,D_2,D_3\}^2}{\sqrt{n}}.
    \]
\end{theorem}

In particular, when $\{X_i\}_{n\in\bN}$ are i.i.d.~the bound becomes $\log n/\sqrt{n}$. In Section~\ref{sec:optimal}, we show that in i.i.d.~case, one can actually remove the $\log n$ factor and conjecture that it is true whenever $\{D_i\}_{i\in\{1,2,3\}}$ are bounded by some constant.

Now, we present the multivariate version involving the joint distribution of sample quantiles. For two vectors $\vx,\vy\in\dR^{\ell}$ we say that $\vx\preceq\vy$ if for all $i\in [\ell]$ we have that $x_i\le y_i.$

\begin{theorem}\label{theorem: multi quantile gen}
	Let $(X_i)_{i\in[n]}$ be a sequence of random variables satisfying Assumption~\ref{as1} and~\ref{as:dn}. Let $A$ be as in \eqref{eq:unifbdd}, $\gTh_n,\theta_n,\mvgm$ be as in~\eqref{def:c}, $D_i$ as in \eqref{eq:D}, $\gS_n$ as in \eqref{def:sigmas}, 
    \begin{equation}\label{eq:mins}
        \theta_{\min}:=\min_{k\in[\ell]}\theta_n(m_k),\quad \gs_{\max}^2:=\max_{k\in[\ell]}\gS_n(k,k) \quad \text{ and } \quad\vQ_n:=\left(Q^{(\ga_1)}_n,Q^{(\ga_2)}_n,\ldots,Q^{(\ga_\ell)}_n\right),
        \end{equation} and  $\mvZ\sim\N(0,I_{\ell}).$
        Assume that there are constants $0<c<C<\infty$ such that for all $n\in\bN$
    \begin{equation}\label{eq:multbdd}
        c<\gs_{\max}^2,\theta_n(m_1),\ldots,\theta_n(m_\ell),\norm{\gS_n^{-\half}}_{\op }<C\quad\text{ and }\quad \max\{D_1,D_2,D_3\}\le C\cdot n^{\frac{1}{4}}.
    \end{equation}
    Then
    \begin{align*}
		 &\sup_{\mvx\in\dR^{\ell}}\left|\pr\left(\Theta_n\sqrt{n}\cdot(\vQ_n-\mvgm)\preceq\mvx\right)-\pr\left(\gS_n^{\half}\mvZ\preceq\mvx\right)\right|\\
		 &\qquad\lesssim D_1\left(\sqrt{\log n+D_1} \right) \cdot\frac{\log n}{\sqrt{n}}
         +\frac{\ell^{1/4}}{{\sqrt{n}}}D_1\left(D_2+D_3\ell^{-1}\right).
    \end{align*}
\end{theorem}

\begin{remark}
    Notice that letting $\ell=1$ and assuming that $\mu=0$ reduces the bound from Theorem~\ref{theorem: multi quantile gen} to the one in Theorem~\ref{theorem: univ DG med}.
\end{remark}

\begin{remark}\label{rem:as}
    One can replace the bounded second derivative assumption in~\eqref{eq:unifbdd} by H\"older continuity assumption for the $F_i'$. But this will give a slower rate of convergence in Theorem~\ref{theorem: univ DG med}.
    One can also relax assumption \eqref{eq:unifbdd_thm} and derive the following bound that holds for $n>N(\eps,\theta_n,\gs_n)$
    \begin{align*}
	 	\dks\left(\frac{\theta_n\sqrt{n}}{\gs_n }M_n\,,\,Z\right)
	 	\lesssim 
	\frac{A}{\theta_n^2}\cdot \frac{D_1\vee \gs_n^2}{\gs_n}\left(1+\frac{\theta_n}{\gs_n}\cdot{\sqrt{\frac{D_1\vee \gs_n^2}{\log n}}} \right) \cdot\frac{\log n}{\sqrt{n}}+\frac{D_1(D_2+D_3)}{\gs_n^3\sqrt{n}}.
\end{align*}
    Notice that in multivariate case $c$ and $C$ depend on $\ell$, again relaxing this assumption one could derive the bound for $n$ large enough 
    \begin{align*}
         &\sup_{\mvx\in\dR^{\ell}}\left|\pr\left(\Theta_n\sqrt{n}\cdot(\vQ_n-\mvgm)\preceq\mvx\right)-\pr\left(\gS_n^{\half}\mvZ\preceq\mvx\right)\right|\\
		 &\qquad\lesssim\frac{A}{\theta_{\min}^2}\cdot\ell\norm{\gS_n^{-\half}}_{\op }\left(D_1\vee\gs^2_{\max}\right)\left(1+\norm{\gS_n^{-\half}}_{\op }\theta_{\min}\cdot{\sqrt{\frac{D_1\vee \gs_{\max}^2}{\log n}}} \right) \cdot\frac{\log n}{\sqrt{n}}\\
        &\qquad\qquad+\frac{\ell^{1/4}}{{\sqrt{n}}}\norm{\gS_n^{-\half}}_{\op }^3D_1\left(D_2+D_3\ell^{-1}\right).
    \end{align*}
    Finally, the constant factor in the bound can also be made explicit. Since we did not optimize over it at in any of the steps in the proof, we decided not to pursue this direction.
\end{remark}

While we present Theorem~\ref{theorem: univ DG med} under the assumption of local dependency, our approach can be applied to other dependency structures. The first step of our argument is applying the classical linearization technique that rewrites $\pr(M_n\le x)$ as the probability that the sum of Bernoulli random variables is less than some value. This step can be done under any dependence of $X_i$'s. However, the mean of these Bernoulli random variables $Y_i$ depends on $x$. While it is not an issue in the derivation of regular CLT, it poses significant complications in bounding the convergence rate when the value of $x$ is ``far'' from the true value of the median. Hence, we consider two cases when $\abs{x}\le K_n$ and $\abs{x}> K_n$ for some threshold function $K_n$. The next step is applying the appropriate approach of Stein's method to treat the former case and establishing a concentration bound for Bernoulli random variables to treat the latter. In the end, we optimize over $K_n$ to derive the result. In conclusion, our argument can be adopted whenever Stein's method can be used to derive CLT for the corresponding Bernoulli random variables $Y_i$, and one can derive a needed concentration inequality.

\subsection{Organization}

This paper is organized as follows. We present an application of our univariate result to the moving-average model in Section~\ref{sec:MA}, and state and prove preliminary results in Section~\ref{sec:prelim}. The proofs of Theorems~\ref{theorem: univ DG med} and ~\ref{theorem: multi quantile gen} are presented in Section~\ref{sec:main. proof}. Finally, we discuss the optimality of the rate of convergence in Section~\ref{sec:optimal}.


\section{Application to the moving-average model}\label{sec:MA}
In this section, we present an application of Theorem~\ref{theorem: univ DG med} to the moving-average model. We refer to~\cite{BrockwellDavis} for an overview of moving-average models and their significance in time series analysis.

One usually considers a sequence $(\gz_i)_{i\in\dZ}$ of i.i.d.~standard normal random variables.  However, to apply our results, we assume that $\{\gz_{i}\}_{i\in\dZ}$ is a sequence of independent (not necessarily identical distributed) random variables that are symmetric around zero and satisfy Assumption~\ref{as1}. In particular, the median of any finite linear combination of $\gz_i$'s is zero.

Thus, one can estimate the median $\mu$ for
the following linear model using sample median based on the observed data $(X_{t})_{t\in[n]}$, where
\begin{align*}
	X_t=\mu+\sum_{i=1}^q c_i\gz_{t-i}+\gz_t
\end{align*}
for all $t\in[n]$. This is also known as $\textrm{MA}(q)$ model with parameters $\mu$, $c_1, c_2,\ldots, c_q$.

Notice that $X_i$ is independent of $X_j$ whenever $\abs{i-j}\ge q+1$. Consider a natural choice of dependency neighborhoods (in this case generated by a graph), where $X_j\in N_i$ if and only if $\abs{i-j}\le q$. Hence the parameters $D_i$ of the dependency neighborhoods, as in \eqref{eq:D}, satisfy $D_1\le 2q+1$, $D_2\le D_1^2$, and $D_3\le D_1^3$. 
Although $X_t$'s are no longer independent, they are still identically distributed. Thus, $\theta$ and $\gs$ are some fixed but unknown constants. In fact, from~\eqref{def:cmed} one can see that $\theta=F_{X_1}'(\mu)$.
Hence, Theorem~\ref{theorem: univ DG med} implies the sample median $M_n(X_1, X_2,\ldots, X_n)$ obeys a central limit theorem with the following upper bound on the rate of convergence,
\begin{align}\label{eq:MA}
	\dks\left(\frac{\theta}{\gs}\cdot \sqrt{n}(M_n-\mu)\,,\,Z\right)\lesssim\frac{q^2+q\log n}{\sqrt{n}},
 \end{align}
where $Z\sim\N(0,1)$.
In particular, the median of $\textrm{MA}(q)$ obeys a Gaussian CLT whenever $q\ll n^{1/4}$. Here, we note that the constant $\theta/\gs$ appearing on the left-hand side of~\eqref{eq:MA} is still an unknown quantity. However, we aim to merely highlight the Gaussian approximation and we leave the estimation of this quantity as a topic for future consideration.


\section{Preliminary Results}\label{sec:prelim}

\subsection{Local dependence}\label{ssec:LD}
In this section, we state the results related to the dependency neighborhoods approach that will be used to prove the main results.

As mentioned above, we will use a CLT for a sum of indicator random variables with local dependencies. Various theorems cover this case in the univariate case (see~\cite{Baldi89, ChenShao04, Rinott94} among others) and give the same order. We chose to present the multivariate version of the result.

\begin{thm}[Theorem~2.1 and Remark 2.2 in~\cite{Fang16}]\label{th: multi DG}
	Suppose that $\mvX_1,\mvX_2,\ldots,\mvX_n$ are $\ell$-dimensional random vectors with $\E \mvX_i =\bf{0}$. Suppose $(\mvX_i)_{i\in[n]}$ have dependency neighborhoods as in Assumption~\ref{as:dn} and for all $i\in[n]$, $j\in N_i$, and $k\in N_{ij}$ we have that
	\[
		\abs{X_i}\le\gb,\quad\abs{N_i}\le D_1,\quad\abs{N_{ij}}\le D_2,\quad\text{and}\quad\abs{N_{ijk}}\le D_3,
	\]
	then letting $\mvW=\frac{1}{\sqrt{n}}\sum_{i=1}^nX_i$ and its variance-covariance matrix given by $\gS$ we have that
	\begin{align*}
		\sup_{A\in\cC}\left|\pr\left(\mvW\in A\right)-\pr\left(\gS^{\half}\mvZ\in A\right)\right|
		\lesssim n\ell^{1/4}\cdot\norm{\gS^{-\half}}^3_{\op }\cdot\gb^3\cdot D_1\left(D_2+D_3/\ell\right),
	\end{align*}
	where $\cC$ is the collection of all convex sets in $\dR^\ell$ and $\norm{\cdot}_{\op}$ denotes the operator norm.
\end{thm}

An important ingredient in our argument is the application of a concentration inequality for locally dependent random variables. Some such inequalities are already present in the literature~\cite{Janson04,paulin2012concentration}. However, they focus on different, often less general, dependency settings. On top of that, a slightly weaker statement than a Hoeffding-type inequality will be sufficient for our purposes.

\begin{lem}\label{th: DN Hoef}
    Let $X_1, X_2,\ldots, X_n$  be a sequence of mean zero random variables taking values in $[-1,1]$. Suppose that $\forall i\in[n]$, $\exists N_i\subseteq[n]$ such that $X_i$ is independent of $(X_l)_{l\in N_i^c}$ (\ie~Assumption~\ref{as:dn}.(1) holds). Then for all $0<t\ll n$, we have
    \[
    \pr\left(\sum_{i=1}^nX_i\ge t\right) \le\exp\left(-\frac{n}{D_1}\cdot g\left(\frac{t}{n}\right)\right),
    \]
    where for $x\in\dR$
    \[
    g(x):=1+x+x^2-e^x.
    \]
\end{lem}
Notice that when $x$ is small $g(x)\approx x^2/2$. Thus, Lemma~\ref{th: DN Hoef} provides a local sub-Gaussian tail bound. 
\begin{proof}
    Consider the moment generating function of $W:=\sum_{i=1}^nX_i$, \ie~  
    $M(r):=\E\left(e^{r W}\right).$ Define $W_{N_i}:=\sum_{j\in N_i}X_j$. Since for every $i\in[n]$ we have that $\E X_i=0$ and $X_i$ is independent of $W-W_{N_i}$ we can rewrite the derivative of $M(r)$ in the following way
    \begin{align*}
        \abs{M'(r)}&\le \sum_{i=1}^n \abs{\E \left(X_i e^{rW}\right)}\\
        &=\sum_{i=1}^n \abs{\E \left(X_i e^{rW} (1-e^{-rW_{N_i}})\right)}
        \le n (e^{rD_1}-1) M(r),
    \end{align*}
 where in the last inequality we used the fact that $\abs{X_i}\le 1$ and $\abs{N_i}\le D_1$ for all $i$. Dividing by $M(r)$ on both sides of the inequality and integrating with respect to $r$ yields that
 \[
 \log M(r)\le n\int_0^r (e^{sD_1}-1)\,ds=\frac{n}{D_1}(e^{rD_1}-1-rD_1).
 \]
This bound allows us to derive the following concentration inequality for $t>0$
\begin{align*}
    \pr\left(W\ge t\right)&\le \exp(-rt)\cdot M(r)
    \le\exp\left(-rt+\frac{n}{D_1}(e^{rD_1}-1-rD_1)\right).
\end{align*}
Picking $r=t/(nD_1)$ and recalling that $t\ll n$ we get the desired inequality
\begin{align*}
    \pr\left(W\ge t\right)
    &\le\exp\left(-\frac{n}{D_1}\left(\frac{t^2}{n^2}+1+\frac{t}{n} -e^{\frac{t}{n}}\right)\right).
\end{align*}
This completes the proof.
\end{proof}

\subsection{Bound on the distance between two Gaussian vectors with the same mean}

One of the terms in the upper bounds in each of our main theorems is of the form
\begin{align}\label{eq:tvbd}
	\sup_{\mvx}\left|\pr(\mvZ_1\preceq\mvx)-\pr(\mvZ_2\preceq\mvx)\right|,
\end{align}
where $\mvZ_1$ and $\mvZ_2$ are $\ell$-dimensional Gaussian random vectors with the same mean and variance--covariance matrices $\gS_1$ and $\gS_2$, respectively. This section is dedicated to addressing this problem.

In the one-dimensional case, one can easily upper bound the Kolmogorov-Smirnov distance by the square root of the $L^{1}$-Wasserstein distance as follows. With $\gr\ge 1$, we have
\begin{align*}
	\dks(\N(0,\gr^2),\N(0,1))
	 &\le\sqrt{4/\sqrt{2\pi}\cdot\dwas(\N(0,\gr^2),\N(0,1))}\\
	 &\le\sqrt{4/\sqrt{2\pi}\cdot\E_{Z'\sim \N(0,\gr^2-1)}\abs{Z'}}
	\le\sqrt{4/\pi}\cdot\abs{\gr^2-1}^{1/4}.
\end{align*}

However, one can derive a better bound by applying Stein's estimates (see Proposition~3.6.1 in~\cite{NourdinPeccati12}) or via direct computation as follows.
\begin{lem}\label{lem:KSbycomp}
	For $\gr\ge 1$, we have
	\begin{align*}
		\dtv(\N(0,\gr^2),\N(0, 1))
		 &\le\frac{2}{\sqrt{\pi}}\abs{\gr-1}\sqrt{\frac{\log\gr}{\gr^2-1}}\cdot\exp\left(-\frac{\log\gr}{\gr^2-1}\right)
		\le\sqrt{\frac{2}{\pi e}}\cdot\abs{\gr-1}.
	\end{align*}
\end{lem}

\begin{proof}
    Letting $\gf$ be the density function of standard normal distribution, we get
	\begin{align*}
		\dtv(\N(0,\gr^2),\N(0,1))
		 & =\int_0^\infty\abs{\gr^{-1}\gf(x/\gr)-\gf(x)}dx
		 = 2\int_{x^*}^{\infty} (\gr^{-1}\gf(x/\gr)-\gf(x)) dx,
	\end{align*}
	where $x^*>0$ satisfies the following equality $\gr^{-1}\exp(-x^2/2\gr^2)=\exp(-x^2/2)$, \ie
	\begin{align}\label{eq:x*}
		x^*=\sqrt{2(\log\gr)/(1-\gr^{-2})}.
	\end{align}
	Letting $Z\sim\N(0,1)$, we get that
	\begin{align*}
		\dtv(\N(0,\gr^2),\N(0,1))
		 & = 2\pr\left(\frac{x^*}{\gr}\le Z\le x^*\right)\\
		 &\le 2x^*\left(1-\gr^{-1}\right)\,\gf\left(x^*/\gr\right) 
		 =\frac{2}{\sqrt{\pi}}(\gr-1)\sqrt{\frac{\log\gr}{\gr^2-1}}\cdot\exp\left(-\frac{\log\gr}{\gr^2-1}\right).
	\end{align*}
	It is easy to check that for all $x\ge 0$ we have $\sqrt{x}e^{-x}\le 1/\sqrt{2e}$. Simplifying, we get the result.
\end{proof}

\begin{rem}
	For general $\gs_1>\gs_2$, letting $\gr=\gs_1/\gs_2> 1$, by the scaling property of normal distribution, we have
	\begin{align*}
		\dtv(\N(0,\gs_1^2),\N(0,\gs_2^2))
		=\dtv(\N(0,\gr^2),\N(0,1)).
	\end{align*}
\end{rem}
In the multivariate case, bounding~\eqref{eq:tvbd} becomes more complex. We refer to~\cite{Devroye18} and references therein for the background on this question. While for our purposes, we need to bound the difference between two measures over specific sets, computing that does not seem feasible. Hence, similarly to the univariate case, we derive a bound in the total variation distance. One can bound this distance in terms of the eigenvalues of the $\gS_1^{-\half}\gS_2\gS_1^{-\half}$ or, equivalently, $\gS_2^{-\half}\gS_1\gS_2^{-\half}$. We present that in Theorem~\ref{theorem:multigaus}, however for clarity of presentation we assume that $\gS_1=I,\gS_{2}=\gS$. Although our technique is simple, to our knowledge, it is not present in the literature, and we find it interesting on its own. For instance, it improves the constant in the upper bound of Theorem 1.1 in~\cite{Devroye18}.
\begin{thm}\label{theorem:multigaus}
    Let $\mvZ_1$ and $\mvZ_2$ be mean zero $\ell$-dimensional normal random vectors with variance-covariance matrices $I$ and $\gS$, respectively. Suppose that $\gl_1,\gl_2,\ldots,\gl_\ell$ are the eigenvalues of $\gS$, then
    $$\dtv(\mvZ_1,\mvZ_2)^2\le\sum_{i=1}^{\ell}\frac{\left(\sqrt{\gl_i}-1\right)^2}{\gl_i+1}\le\norm{\gS-I}_{\F}^2.$$
\end{thm}

The following bound will be helpful when none of the variance--covariance matrices are identity matrices.

\begin{lem}\label{lem:opbdd}
	Let $\gS_1$ be a $\ell\times \ell$ invertible symmetric matrix and $\gS_2$ be another $\ell\times \ell$ matrix. Then
	\[
		\norm{\gS_1^{-\half}\gS_2\gS_1^{-\half}}_{\F}\le\norm{\gS_1^{-1}}_{\op}\cdot\norm{\gS_2}_{\F}.
	\]
\end{lem}
\begin{proof}
	Let $\gl_1,\gl_2,\ldots,\gl_{\ell}$ be the eigenvalues of $\gS_1$, and $\{u_i\mid 1\le i\le \ell\}$ be the corresponding orthonormal eigenbasis. Then $\gS_1^{-1}$ can be written as
	$
		\gS_1^{-1}=\sum_{i=1}^d\gl_i^{-1}u_iu_i^{\T}.
	$
	Using this, we get that
	\begin{align*}
		\norm{\gS_1^{-\half}\gS_2\gS_1^{-\half}}_{\F}^2
		 & =\sum_{i,j=1}^{\ell}(\gl_i\gl_j)^{-1}(u_j^{\T}\gS_2u_i)^2
		 \le\norm{\gS_1^{-1}}_{\op}^2\cdot\sum_{i,j=1}^{\ell}(u_j^{\T}\gS_2u_i)^2
		=\norm{\gS_1^{-1}}_{\op}^2\cdot\norm{\gS_2}^2_{\F}.
	\end{align*}
	This completes the proof.
\end{proof}

Let $f,g$ be the density of the random variables $X, Y$, respectively, and recall the definitions of the Hellinger distance and affinity (also referred to as the Bhattacharyya coefficient) defined as
    $$\dhell(X,Y):=\left(\frac12\int \left(\sqrt{f}-\sqrt{g}\right)^2\right)^{\frac12} \quad \text{ and } \quad \ga(X,Y):=\int\sqrt{fg}.$$
    It is well known~\cite{Strasser1985} that the total variation distance can be upper bounded in terms of these distances as follows
    \begin{equation}\label{eq:distancesOri}
        \dtv(X,Y)\le \sqrt{2}\dhell(X,Y)=\sqrt{2}\cdot\sqrt{1-\ga(X,Y)}.
    \end{equation}
    In fact, a simple computation yields a slightly stronger inequality, namely
     \begin{equation}\label{eq:distances}
        \dtv(X,Y)\le \sqrt{1-\ga(X,Y)^2}.
    \end{equation}
    Indeed,
    \begin{align*}
        \dtv(X,Y)^2&=\left(\int f\vee g-1\right)\left(1-\int f\wedge g\right)=1-\int f\vee g \cdot \int f\wedge g\\
        &\le1-\left(\sqrt{(f\vee g)(f \wedge g)}\right)^2=1-\alpha(X,Y)^2.
    \end{align*}

\begin{proof}[Proof of Theorem~\ref{theorem:multigaus}]
    
    It is enough to bound $\int\sqrt{fg}$, where
    \begin{align*}
		f(\mvx) & =(2\pi)^{-\ell/2}\cdot\exp\left(-\mvx^{\T}\mvx/2\right)
		\text{ and }
		g(\mvx)=(2\pi)^{-\ell/2}\det(\gS)^{-1/2}\cdot\exp\left(-\mvx^{\T}\gS^{-1}\mvx/2\right).
    \end{align*}
    Then
    \begin{align*}
		\sqrt{f(\mvx)g(\mvx)} & =(2\pi)^{-\ell/2}\cdot\det(\gS)^{-1/4}\cdot\exp\left(-\mvx^{\T}(I+\gS^{-1})\mvx/4\right)\\
		           & =\frac{\det(\gS)^{1/4}}{\det((I+\gS)/2)^{1/2}}\cdot (2\pi)^{-\ell/2}\cdot\frac{\det((I+\gS)/2)^{1/2}}{\det(\gS)^{1/2}}\cdot\exp\left(-\mvx^{\T}(I+\gS^{-1})\mvx/4\right).
    \end{align*}
    Recognizing the density function of a normal distribution with mean zero and variance-covariance matrix $(I+\gS^{-1})/{2}$, we conclude that
    \begin{align*}
		\ga(\mvZ_1,\mvZ_2)=\int\sqrt{fg} & =\left(\frac{\det(\gS)^{1/2}}{\det((I+\gS)/2)}\right)^{\frac12}=\left(\prod_{i=1}^{\ell}\frac{2\sqrt{\gl_i}}{1+\gl_i}\right)^{\frac12}.
    \end{align*}
    Hence, by \eqref{eq:distances} we get that
    \begin{align*}
		\dtv(\mvZ_1,\mvZ_2)^2
		\le 1-\left(\prod_{i=1}^{\ell}\frac{2\sqrt{\gl_i}}{1+\gl_i}\right)
		\le\sum_{i=1}^{\ell}\left(1-\frac{2\sqrt{\gl_i}}{1+\gl_i}\right)
		=\sum_{i=1}^{\ell}\frac{\left(\sqrt{\gl_i}-1\right)^2}{\gl_i+1}.
    \end{align*}
    The conclusion follows from the fact that $\norm{\gS-I}_{\F}^{2}=\sum_{i=1}^{\ell}(\gl_i-1)^2$.
\end{proof}
\subsection{Variance Control}\label{sec:Variance}
We present our main results with the implicit variance $\gs=\gs_{n,0}$. To compute $\gs_{n,0}$ explicitly, one would have to rely on particular properties of the model. For example, as we mentioned before, suppose in addition to Assumptions~\ref{as1}, we have that
$	\pr(X_i\vee X_j\le 0)=1/4$ for all $i\neq j,
$
then $\gs\equiv\frac12$.

\begin{lem}\label{lem:gSgap}
	In the situation of Theorem~\ref{theorem: multi quantile gen}, let $\theta_n$ be as in~\eqref{def:cmed}, and $\eps,A$ are as in~Assumption~\ref{as1}.\ref{as:unifbdd}. Then for any $\mvx\in\dR^{\ell}$ such that $\max_{k\in [\ell]}\abs{x_k/\theta_n(m_k)}\le \eps \sqrt{n}$ we have that
	\begin{align*}
		\dtv\left(\vN_{\ell}\left(0,\gS_{n,\mvx}\right),\,\vN_{\ell}\left(0,\gS_0\right)\right)
		\le \frac{3\sqrt{2}A\ell}{\theta_{\min}}\cdot\norm{\gS_0^{-1}}_{\op}\cdot\frac{D_1\norm{\mvx}_\infty}{\sqrt{n}}.
	\end{align*}
\end{lem}

\begin{proof}
    Applying Theorem~\ref{theorem:multigaus} and Lemma~\ref{lem:opbdd} gives us that
	\begin{align*}
		\dtv\left(\vN_{\ell}\left(0,\gS_{n,\mvx}\right),\,\vN_{\ell}\left(0,\gS_0\right)\right)
		 &\le\norm{\gS_0^{-\half}{\gs_n}_{n,\mvx}\gS_0^{-\half}-I}_{\F}\\
		 &\le\norm{\gS_0^{-1}}_{\op}\cdot\norm{\gS_{n,\mvx}-\gS_0}_{\F}.
	\end{align*}
	To compute $\norm{\gS_{n,\mvx}-\gS_0}_{\F}$, we bound the difference in each coordinate. Define the function
	$$
		f_{k,\mvx}(u):=\ind_{u\in(m_{k}, m_{k}+n^{-\half}\cdot x_k/\theta_n(m_k))}-\ind_{u\in(m_{k}+n^{-\half}\cdot x_k/\theta_n(m_k),m_{k})}
	$$
	for $k\in[\ell]$. Then
	\[
		\1_{u\le m_k+n^{-\half}\cdot x_k/\theta_n(m_k)} =\1_{u\le m_k} +f_{k,\mvx}(u).
	\]
    In particular, using Assumption~\ref{as1}.\ref{as:unifbdd} together with MVT we derive that
    \begin{align}\label{eq:Ef}
    \E \abs{f_{k,\mvx}(X_j)} &\le \pr\left(X_j\in (m_{k}, m_{t}+n^{-\half}\cdot x_k/\theta_n(m_k))\right)+\pr\left(X_j\in(m_{k}+n^{-\half}\cdot x_k/\theta_n(m_k),m_{k})\right)\notag\\
    &\le A\cdot\frac{D_1}{\sqrt{n}}\abs{\frac{x_k}{\theta_n(m_k)}}.
    \end{align}
	Thus for $s$ and $t$ in $[\ell]$,  we get 
	\begin{align*}
		 & \abs{(\gS_{n,\mvx})_{st}-\gS_{st}}\\
		 & =\frac{1}{n}\abs{\sum_{i,\,j\in N_{i}}\cov\left(\1_{X_i\le m_{s}}, f_{t,\mvx}(X_j)\right)+\cov\left(f_{s,\mvx}(X_i),\1_{X_t\le m_t}\,\right)+\cov\left(f_{s,\mvx}(X_i), f_{t,\mvx}(X_j)\right)}\\
		 &\le 2A\cdot\frac{D_1}{\sqrt{n}}\cdot\left(\abs{\frac{x_s}{\theta_s}}+\abs{\frac{x_t}{\theta_t}}+\sqrt{\abs{\frac{x_sx_t}{\theta_s\theta_t}}}\right)
		\le3A\cdot\frac{D_1}{\sqrt{n}}\cdot\left(\abs{\frac{x_s}{\theta_s}}+\abs{\frac{x_t}{\theta_t}}\right),
	\end{align*}
    where in the first inequality, we bounded each indicator $\1_{X_i\le m_{k}}$ by $1$, used \eqref{eq:Ef}, and applied the Cauchy–Schwarz inequality for the third term within the sum.
	Hence
	\begin{equation}\label{eq:var comp}
		\norm{\gS_{n,\mvx}-\gS_0}_{\F}\le3A\cdot\frac{D_1}{\sqrt{n}}\cdot\sqrt{\sum_{s,t}\left(\abs{\frac{x_s}{\theta_s}}+\abs{\frac{x_t}{\theta_t}}\right)^2}\le 3\sqrt{2}\ell\cdot\frac{A}{\theta_{\min}}\cdot\frac{D_1 }{\sqrt{n}}\cdot\norm{\mvx}_{\infty}.
	\end{equation}
	This completes the proof.
\end{proof}
\begin{rem}
	In Lemma~\ref{lem:gSgap}, $D_1$ can be replaced by the average size of the neighborhood instead of the maximum one.
\end{rem}


\section{Proof of Main Results}\label{sec:main. proof}

\subsection{Limiting distribution of the sample median}\label{sec:univ. proof}

\begin{proof}[Proof of Theorem~\ref{theorem: univ DG med}]

 Recall the definition of the sample median
	$$
	M_n=\sup\left\{x\in\dR:\abs{\{i\mid X_i\le x\}}\le\lfloor n/2\rfloor\right\}
	$$
 and so
 \begin{align*}
     \pr\left(M_n\le  y\right)&=\pr\left(\abs{\{i\mid X_i \le y\}} \ge \lfloor{n}/{2}\rfloor+1\right)
		 =\begin{cases}
		     \pr\left(\abs{\{i\mid X_i\le y\}}\ge \frac{n+2}{2}\right) &\textnormal{if $n$ is even},\\
            \pr\left(\abs{\{i\mid X_i\le y\}}\ge \frac{n+1}{2}\right) &\textnormal{if $n$ is odd}.
		 \end{cases}
 \end{align*}
 Thus fixing $x\in\dR$. We have
	\begin{align}
		\pr\left(\theta_n\sqrt{n}M_n\le x\right)
		 &{=}\pr\left(M_n\le n^{-\half}\cdot x/\theta_n\right)\notag\\
		 &\in\left[\pr\left(\abs{\{i\mid X_i\le n^{-\half}\cdot x/\theta_n\}}\ge {n}/{2}+1\right),\pr\left(\abs{\{i\mid X_i\le n^{-\half}\cdot x/\theta_n\}}\ge {n}/{2}\right)\right]\notag\\
		 & =\left[\pr\left(\sum_{i=1}^n\1_{X_i\le n^{-\half}\cdot x/\theta_n}\ge {n}/{2}+1\right),\pr\left(\sum_{i=1}^n\1_{X_i\le n^{-\half}\cdot x/\theta_n}\ge {n}/{2}\right)\right]\label{eq:step 1}
	\end{align}
	By subtracting $\sum_{i=1}^nF_i\left(n^{-\half}\cdot x/\theta_n\right)$ to each side of the inequality inside of the probability 
	we get
	\begin{align}\label{eq:step 2}
		 &\pr\left(\sum_{i=1}^n\1_{X_i\le n^{-\half}\cdot x/\theta_n}\ge {n}/{2}\right)
		 =\pr\left(\sum_{i=1}^n\left(\1_{X_i\le n^{-\half}\cdot x/\theta_n}-F_i(n^{-\half}\cdot x/\theta_n)\right)\ge -\sqrt{n}\cdot x_n\right)
	\end{align}
	where
	\[
	x_n = n^{-\half}\cdot\left(\sum_{i=1}^nF_i(n^{-\half}\cdot x/\theta_n) -{n}/{2}\right).
	\]
	By Assumption~\ref{as1}, $F_i$'s are twice continuously differentiable at $0$. Thus, using Taylor expansion for each $F_i(x)$ at $0$, $x_n$ could be rewritten as 
	\begin{align*}
		x_n=\frac{1}{\sqrt{n}}\left({n}/{2}-{n}/{2}\right)+\frac{x}{\theta_n n}\sum_{i=1}^n F_i'(0)+\frac{x^2}{2\theta_n^2 n^{3/2}}\sum_{i=1}^nF_i''(y_i),
	\end{align*}
    where $y_i\in (0,x/(\theta_n \sqrt{n}))$. Recall constants $\eps$, $A$ from Assumption~\ref{as1}.\ref{as:unifbdd} and constants $c,C$ from~\eqref{eq:unifbdd_thm}. 
    
    Define 
    \[
    K_n^2=\frac32\cdot \frac{\gs_n^2}{\min(1,c)^2}\cdot  D_1\cdot \log n,  
    \]
    so that 
    \begin{align}\label{eq:Kc}
        K_n^2/(3D_1)\ge \frac12\log n,\quad K_n^2/(2\gs_n^2)\ge \frac12\log n
    \end{align}
    and 
    \[
    K_n^2\le \frac{3C^3}{2c^2}\cdot n^{\frac14}\log n \le \frac{3C^3}{2c^3} \cdot \theta_n\cdot n^{\frac14}\log n\le \eps \cdot \theta_n\cdot \sqrt{n}
    \]
    holds for 
    \begin{equation}\label{eq:N(eps,c,C)}
    n\ge N(\eps,c,C), \textnormal{ where } N(\eps,c,C) \textnormal{ is the smallest integer solution to } \frac{n^{1/4}}{\log n}\ge \frac{3C^3}{2\eps c^3}.
    \end{equation}
    Thus, when $\abs{x}\le K_n$ we have $\max_{i\in[n]}\abs{y_i}\le \eps$ and from Assumption~\ref{as1}.\ref{as:unifbdd} it follows that for all $n\ge N(\eps,c,C)$
	\begin{align}\label{eq:TaylorR}
		\abs{x_n-x}\le\frac{K_n^2}{2\sqrt{n}}\cdot\frac{A}{\theta_n^2}.
	\end{align}
    
	On the other hand, the left-hand side of the inequality inside of probability in~\eqref{eq:step 2} is a scaled sum of centered Bernoulli random variables $$
	Y_{i,x}:=\1_{X_i\le n^{-\half}\cdot x/\theta_n}-F_i(n^{-\half}\cdot x/\theta_n),\quad i\in[n].
	$$
	
	We treat the left bound of the interval from~\eqref{eq:step 1} similarly. This allows us to rewrite the original quantity of interest in the following way.
	\begin{align}\label{eq:linear}
		\pr\left(\theta_n\sqrt{n}M_n\le x\right)\in\left[\pr\left(\sum_{i=1}^nY_{i,x}\ge 1-\sqrt{n}\cdot x_n\right),\pr\left(\sum_{i=1}^nY_{i,x}\ge -\sqrt{n}\cdot x_n\right)\right].
	\end{align}
	Moreover, for each $x\in\dR$ random variables $\{Y_{i,x}\}$ have the same dependency structure as $\{X_i\}$. 
    We will now consider two cases when $|x|> K_n$, in which we apply a concentration inequality, and $|x|\le K_n$, where we use Stein's method for normal approximation.

The Kolmogorov--Smirnov distance between $({\theta_n\sqrt{n}}/{\gs_n})M_n$ and the standard normal random variable $Z$ can be approximated in terms of $\{Y_{i,x}\}_{i\le n}$ as follows
	\begin{align}\label{eq:max form}
	\begin{split}
	  	&\dks\left(\frac{\theta_n\sqrt{n}}{\gs_n}\cdot M_n\,,\,Z\right)=\sup_{x\in\dR}\left|\pr\left(\theta_n\sqrt{n}\cdot M_n\le x\right)-\Phi(x/\gs_n)\right|\\
		&\le\sup_{x\in\dR}\left\{\left|\pr\left(\sum_{i=1}^nY_{i,x}\ge 1-\sqrt{n}\cdot x_n\right)-\Phi(x/\gs_n)\right|,\left|\pr\left(\sum_{i=1}^nY_{i,x}\ge -\sqrt{n}\cdot x_n\right)-\Phi(x/\gs_n)\right|\right\}.
\end{split}
	\end{align}
	Focusing on the second term inside of the maximum in~\eqref{eq:max form} we rewrite it as
	\begin{align*}
		\sup_{x\in\dR}\left|\pr\left(\sum_{i=1}^nY_{i,x}\ge -\sqrt{n}\cdot x_n\right)-\Phi(x/\gs_n)\right| &\leq\err_1+\err_2,
	\end{align*}
	where
	\begin{align}
		\err_1
		& =\sup_{\abs{x}\ge K_n}\left|\pr\left(\sum_{i=1}^nY_{i,x}\ge -\sqrt{n}\cdot x_n\right)-\Phi(x/\gs_n)\right|\label{eq:E1},\\
		\text{and}\quad\err_2 & =\sup_{\abs{x}<K_n}\left|\pr\left(\sum_{i=1}^nY_{i,x}\ge -\sqrt{n}\cdot x_n\right)-\Phi(x/\gs_n)\right|\label{eq:E2}.
	\end{align}

	\noindent\textit{Bound on $\err_1$}:
	We bound the first term by considering two cases: when $x_n$ is negative and when it is positive. When $x_n<0$ both quantities in~\eqref{eq:E1} are negligible. Indeed, by Lemma~\ref{th: DN Hoef} we have the following concentration bound
	\begin{align*}
		\pr\left(\sum_{i=1}^nY_{i,x}\ge -\sqrt{n}\cdot x_n\right)\le\pr\left(\sum_{i=1}^nY_{i,x}\ge \sqrt{n}\cdot K_n\right)\lesssim\exp\left(-{K_n^2}/{3D_1}\right) {\le n^{-1/2}}
	\end{align*}
	and since $\max\left\{\Phi_{\gs_n^2}(-K_n),1-\Phi_{\gs_n^2} (K_n)\right\}\le\exp\left(-K_n^2/(2\gs_n^2)\right) {\le n^{-1/2}}$ we get that
	\begin{align*}
		\sup_{x<-K_n}\left|\pr\left(\sum_{i=1}^nY_{i,x}\ge -\sqrt{n}\cdot x_n\right)-\Phi(x/\gs_n)\right|\le  {n^{-1/2}}
	\end{align*}
	On the other hand, when $x_n>0$, both terms are close to $1$,
	\begin{align*}
		\sup_{x>K_n}&\abs{\pr\left(\sum_{i=1}^nY_{i,x}\ge -\sqrt{n}\cdot x_n\right)-\Phi(x/\gs_n)}\\
        & =\sup_{x>K_n}\pr\left(\sum_{i=1}^nY_{i,x}\le -\sqrt{n}\cdot x_n\right)+\Phi(-x/\gs_n)\lesssim  {n^{-1/2}}.
	\end{align*}
	
	\noindent\textit{Bound on $\err_2$}:
	 We now consider the second error term with $|x|\le K_n$. Recall that $$\gs_{n,x}^2=\frac1n\var\left(\sum_{i=1}^nY_{i,x}\right).$$ 
     Then $\err_2$ can be bounded as follows

	\begin{align}
		\err_2 &\leq\sup_{\abs{x}\le K_n}\left|\pr\left(\sum_{i=1}^nY_{i,x}\ge -\sqrt{n}\cdot x_n\right)-\Phi(x_n/\gs_{n,x})\right|\label{term: CLT}\\
		    &\qquad\qquad\qquad+\sup_{\abs{x}\le K_n}\left|\Phi(x_n/\gs_{n,x})-\Phi(x_n/\gs_n)\right|\label{term: gaus diff}\\
		    &\qquad\qquad\qquad\qquad\qquad+\sup_{\abs{x}\le K_n}\left|\Phi(x_n/\gs_n)-\Phi(x/\gs_n)\right|.\label{term: gaus x}
	\end{align}

    The terms \eqref{term: gaus diff} and \eqref{term: gaus x} are easily bound using the previously established inequalities. First, we established that 
	\begin{equation}\label{eq:gausdiff}
	~\eqref{term: gaus diff}\le\sup_{\abs{x}\le K_n}\frac{\abs{\gs_n^2-\gs_{n,x}^2}}{\gs_n^2-\max\{0,\gs_n^2-\gs_{n,x}^2\}}\le\sup_{\abs{x}\le K_n}\frac{6\sqrt{2}A}{\theta_n}\cdot\frac{D_1\abs{x}}{\gs_n^2\sqrt{n}}\lesssim\frac{D_1\, K_n}{\sqrt{n}} {\lesssim D_1^{\frac32}\cdot\sqrt{\frac{\log n}{n}}},
	\end{equation}
        where the first inequality follows from Lemma~\ref{lem:KSbycomp} and the second inequality from the variance comparison as in \eqref{eq:var comp}.
        
	Turning our attention to the term~\eqref{term: gaus x}, we notice that as long as $K_n\le\eps \cdot\theta_n \sqrt{n}$ we can use the bound from \eqref{eq:TaylorR} to derive that
	\begin{align}\label{eq:gaus x}
		~\eqref{term: gaus x}\le\sup_{\abs{x}\le K_n}\frac{\abs{x_n-x}}{\gs_n}\le\frac{A}{\theta_n^2}\cdot\frac{K_n^2}{\gs_n\sqrt{n}} {\lesssim D_1\cdot\frac{\log n}{\sqrt{n}}}.
	\end{align}

	Finally, it remains to bound \eqref{term: CLT}. Since all random variables $|Y_i|\le 1$, by Theorem~\ref{th: multi DG} we have a quantitative version of CLT for $\{Y_i\}_{i\ge 1}$. Since $Y_{i,x}$ are centered we we can rewrite~\eqref{term: CLT} as
	\begin{align*}
		\sup_{y\in\dR}\abs{\pr\left(\sum_{i=1}^nY_{i,x}\le\sqrt{n}\cdot y\right)-\Phi(y/\gs_{n,x})}
		 &\lesssim\frac{D_1(D_2+D_3)}{\gs_{n,x}^3\sqrt{n}}\lesssim\frac{D_1(D_2+D_3)}{\gs_n^3\sqrt{n}}\lesssim \frac{D_1(D_2+D_3)}{\sqrt{n}},
	\end{align*}
        where we again used the variance comparison \eqref{eq:var comp} to justify the change from $\gs_{n,x}$ to $\gs_n$ for all $x$ such that $\abs{x}\le K_n \le \eps \cdot\theta_n \sqrt{n}$ and then \eqref{eq:unifbdd_thm} to bound $\gs_n$ by a universal constant.
		Putting these bounds together, we conclude that 
        \begin{align*}
            \sup_{x\in\dR} &\left|\pr\left(\sum_{i=1}^nY_i\ge -\sqrt{n}\cdot x_n\right)-\Phi(x/\gs_n)\right|\notag\\
		&\lesssim n^{-1/2}+D_1^{\frac32}\cdot\sqrt{\frac{\log n}{n}}+D_1\frac{\log n}{\sqrt{n}}+\frac{D_1(D_2+D_3)}{\sqrt{n}}\\
        &\lesssim D_1\sqrt{\log n + D_1} \cdot\sqrt{\frac{\log n}{n}}+\frac{\max\{D_1,D_2,D_3\}^2}{\sqrt{n}},
        \end{align*}
        for $n\ge N(\eps,c,C)$. Multiplying the right-hand side by a universal constant that is large enough yields the bound for all $n$.
	\end{proof}

 \subsection{Limiting joint distribution of sample quantiles}

\begin{proof}[Proof of Theorem~\ref{theorem: multi quantile gen}]\label{sec:ApB}
Define $$\overline{\vQ}_n:=\vQ_n-\mvgm.$$

	Since each centered empirical quantile can be rewritten as
	$$
	\overline{Q}_{n}^{(\ga_k)}:=\sup\left\{x\in\dR:\abs{\{i:X_i\le x\}}\le\left\lfloor n\ga_k\right\rfloor\right\}-m_{k},
	$$
	for any $\mvx=(x_1,x_2,\ldots,x_\ell)\in\dR^\ell$ we rewrite the joint distribution of quantiles in the following way
	\begin{align}
	\begin{split}\label{eq:step 1a}
	  	\pr(\Theta_n\sqrt{n}&\cdot(\vQ_n-\mvgm)\preceq\mvx)\\
		  & =\pr\left(\theta_1\sqrt{n}\,\overline{Q}^{(\ga_1)}_n\le x_1,\ldots,\theta_\ell\sqrt{n}\,\overline{Q}^{(\ga_\ell)}_n\le x_\ell\right)\\
		  &\in\Biggl[\pr\left(\sum_{i=1}^n\1_{X_i-m_k\le (\theta_n(m_k)\sqrt{n})^{-1}x_k}\ge\left\lfloor n\ga_k\right\rfloor+1,\textrm{ for all } k\in[\ell]\right),\\
		  &\qquad\qquad\pr\left(\sum_{i=1}^n\1_{X_i-m_k\le (\theta_n(m_k)\sqrt{n})^{-1}x_k}\ge\left\lfloor n\ga_k\right\rfloor,\textrm{ for all } k\in[\ell]\right)\Biggr].
	\end{split}
	\end{align}
	Following the same idea as in the proof of Theorem~\ref{theorem: univ DG med} we treat each of the probabilities from~\eqref{eq:step 1a} similarly but separately. In particular, by subtraction corresponding values of CDF functions $F_i$ on both sides of the inequalities inside of the probabilities and then using Taylor expansion, we get an analogs equation to~\eqref{eq:linear}, namely
	\begin{align}\label{eq:cdf}
		\pr\left(\Theta_n\cdot (\vQ_n-\mvgm)\preceq\mvx\right) &\in\left[\pr\left(\sum_{i=1}^n\mvY_{i,\mvx}\succeq \mv1-\sqrt{n}\cdot\mvx_n\right),\pr\left(\sum_{i=1}^n\mvY_{i,\mvx}\succeq-\sqrt{n}\cdot\mvx_n\right)\right],
	\end{align}
	where $\mv1:=(1,1,\ldots,1)$ and
	\begin{align*}
		\mvx_n:=\frac{1}{\sqrt{n}}\left(\sum_{i=1}^nF_i\left(m_k+n^{-\half}\cdot x_k/\theta_n(m_k)\right) -\left\lfloor n\ga_k\right\rfloor\right)_{k=1}^\ell.
	\end{align*}
    Define $K_n$ in a similar fashion to the univariate case replacing $\gs_n$ by $\gs_{\max}$
    \[
    K_n^2=\frac32\cdot \frac{\gs_{\max}^2}{\min(1,c)^2}\cdot  D_1\cdot \log n,  
    \]
    then if $\norm{\mvx}_\infty\le K_n$ we have that $\mvx_n-\mvx$ satisfies the following inequality for all for all $n\ge N(\eps,c,C)$, defined in \eqref{eq:N(eps,c,C)},
	\begin{align}\label{eq:xn error}
		\norm{\mvx_n-\mvx}_2\le\frac{\norm{\mvx}_2^2}{2\sqrt{n}}\cdot\frac{A}{\theta_{\min}^2}
		\le\frac{\ell\norm{\mvx}_\infty^2}{2\sqrt{n}}\cdot\frac{A}{\theta_{\min}^2}.
	\end{align}

	For some function $K_n$, each of the terms on the right-hand side of~\eqref{eq:cdf} are treated similarly. For simplicity, we present the bound on the right bound of the interval.
	\begin{align*}
		\sup_{\mvx} &\abs{\pr\left(\sum_{i=1}^n\mvY_{i,\mvx}\succeq-\sqrt{n}\cdot\mvx_n\right)-\pr\left(\gS_n^{\half}\mvZ\succeq-\mvx\right)}\le\wh{\err}_1+\wh{\err}_2,
	\end{align*}
	where
	\begin{align}
		\wh{\err}_1         & =\sup_{\norm{\mvx}_{\infty}\ge K_n}\abs{\pr\left(\sum_{i=1}^n\mvY_{i,\mvx}\succeq-\sqrt{n}\cdot\mvx_n\right)-\pr\left(\gS_n^{\half}\mvZ\succeq-\mvx\right)}\label{eq:multiErr_1}\\
		\text{and}\quad\wh{\err}_2 & =\sup_{\norm{\mvx}_{\infty}< K_n}\abs{\pr\left(\sum_{i=1}^n\mvY_{i,\mvx}\succeq-\sqrt{n}\cdot\mvx_n\right)-\pr\left(\gS_n^{\half}\mvZ\succeq-\mvx\right)}\label{eq:multiErr_2}
	\end{align}
	\noindent\textit{Bound on $\wh{\err}_1$.} Recall that $\gs_{\max}^2:=\max_{i\in[\ell]}\gS_n(i,i)$.
 
	In the first term, we bound similarly to the way we bounded~\eqref{eq:E1} in the one-dimensional case. If $x_{\min}:=\min_k{x_k}\le -K_n$ then using a union bound and Lemma~\ref{th: DN Hoef},
	\begin{align*}
		\sup_{x_{\min}\le -K_n}\abs{\pr\left(\sum_{i=1}^n\mvY_{i,\mvx}\succeq-\sqrt{n}\cdot\mvx_n\right)-\pr\left(\gS_n^{\half}\mvZ\succeq-\mvx\right)}\le 2\ell\exp\left({-{K_n^2}/{(3D_1\vee 2\gs_{\max}^2)}}\right).
	\end{align*}

	On the other hand if $x_{\min}> -K_n$ and $x_{\max}:=\max_k{x_k}\ge K_n$, then we divide coordinates into two sets $A:=\{k:x_k\in(-K_n,K_n)\}$ and  $B:=\{k: x_k\le -K_n\}$. For simplicity let
	\[
		E_k:=\left\{\sum_{i=1}^n(\mvY_{i,\mvx})_k\ge -\sqrt{n}\cdot(\mvx_n)_k\right\}\quad\text{and}\quad F_k:=\left\{\gS_n^{\half}Z_k\ge-(x)_k\right\}.
	\]

	Thus, for all $\mvx$ such that $x_{\min}> -K_n$ and $ x_{\max}\ge K_n$ we have

	\begin{align}
		 &\abs{\pr\left(\sum_{i=1}^n\mvY_{i,\mvx}\succeq-\sqrt{n}\cdot\mvx_n\right)-\pr\left(\gS_n^{\half}\mvZ\succeq-\mvx\right)}\notag\\
		 &\qquad=\abs{\pr\left(\bigcap_{k\in A}E_k\bigcap_{k'\in B}E_{k'}\right)-\pr\left(\bigcap_{k\in A}F_k\bigcap_{k'\in B} F_{k'}\right)}\notag\\
		 &\qquad\le\abs{\pr\left(\bigcap_{k=1}^{\ell}E_k\right)-\pr\left(\bigcap_{k=1}^{\ell}F_k\right)}\label{eq:tailterm1}\\
		 &\qquad\quad+\abs{\pr\left(\bigcap_{k\in A}E_k\left(\bigcap_{k'\in B}E_{k'}\right)^{c\,}\right)-\pr\left(\bigcap_{k\in A} F_k\left(\bigcap_{k'\in B} F_{k'}\right)^{c\,}\right)}\label{eq:tailterm2},
	\end{align}
	where in~\eqref{eq:tailterm2}, we switched to the compliment events. This enables the application of a concentration inequality
	\begin{align*}
		~\eqref{eq:tailterm2} &\le\abs{\pr\left(\bigcup_{k\in B}\left\{\sum_{i=1}^n\mvY_{i,x_k}< -\sqrt{n}\cdot(x_n)_k\right\}\right)-\pr\left(\bigcup_{k\in B}\left\{\gS_n^{\half}Z_k<-(x)_k\right\}\right)}\\
		&\lesssim 2\ell\exp\left({-K_n^2/(3D_1\vee 2\gs^2_{\max})}\right)\lesssim \ell n^{-1/2}.
	\end{align*}
 
	Finally, the remaining term~\eqref{eq:tailterm1} can be bounded similarly to the term $\wh{\err}_2$ as presented below.\\

        \noindent\textit{Bound on $\wh{\err}_2$.} We rewrite $\wh{\err}_2$ in a similar way to how we treated $\err_2$ in the univariate case
	\begin{align}
		\wh{\err}_2
		 & =\sup_{\norm{\mvx}_{\infty}\le K_n}\abs{\pr\left(\sum_{i=1}^n\mvY_{i,\mvx}\succeq -\sqrt{n}\cdot\mvx_n\right)-\pr\left(\gS_{n,\mvx}^{\half}\mvZ\succeq-\mvx_n\right)}\label{term: multiCLT}\\
		 &\quad+\sup_{\norm{\mvx}_{\infty}\le K_n}\abs{\pr\left(\gS_{n,\mvx}^{\half}\mvZ\succeq-\mvx_n\right)-\pr\left(\gS_n^{\half}\mvZ\succeq-\mvx_n\right)}\label{term: multigaus diff}\\
		 &\quad+\sup_{\norm{\mvx}_{\infty}\le K_n}\abs{\pr\left(\gS_n^{\half}\mvZ\succeq-\mvx_n\right)-\pr\left(\gS_n^{\half}\mvZ\succeq-\mvx\right)}.\label{term: multigaus x}
	\end{align}
	Applying multivariate CLT as in Theorem~\ref{th: multi DG} we get
	\begin{align*}
		~\eqref{term: multiCLT} &\le\sup_{\norm{\mvx}_{\infty}\le K_n}\sup_{A\in\cC}\abs{\pr\left(\frac{1}{\sqrt{n}}\sum_{i=1}^n\mvY_{i,\mvx}\in A\right)-\pr\left(\gS_{n,\mvx}^{\half}\mvZ\in A\right)}\\
		&\lesssim\frac{\ell^{1/4}}{{\sqrt{n}}}\norm{\gS_{n,\mvx}^{-\half}}^3_{\op }D_1\left(D_2+D_3\ell^{-1}\right)
        \lesssim\frac{\ell^{1/4}}{{\sqrt{n}}}\norm{\gS_n^{-\half}}^3_{\op }D_1\left(D_2+D_3\ell^{-1}\right)\\
        &\lesssim\frac{\ell^{1/4}}{{\sqrt{n}}}D_1\left(D_2+D_3\ell^{-1}\right)
	\end{align*}

	To bound the error term~\eqref{term: multigaus diff} we derive the bound on the total variation distance and invoke Lemma~\ref{lem:gSgap}, which states that
	\begin{align*}
		~\eqref{term: multigaus diff}
        &\le\sup_{\norm{\mvx}_{\infty}\le K_n}\max\{\norm{\gS_n^{-1}}_{\op},\norm{\gS_{n,\mvx}^{-1}}_{\op}\}\cdot \norm{\gS_{n,\mvx}-\gS_n}_{\F}\notag\\
		&\le\sup_{\norm{\mvx}_{\infty}\le K_n}\frac{6\sqrt{2}A}{\theta_{\min}}\cdot\ell\cdot\norm{\gS_n^{-1}}_{\op}\cdot\frac{D_1\norm{\mvx}_\infty}{\sqrt{n}}
		\lesssim\frac{D_1K_n}{\sqrt{n}}\lesssim D_1^{\frac32}\sqrt{\frac{\log n}{n}}.\notag
	\end{align*}

	Finally, the term~\eqref{term: multigaus x} we bound coordinate-wise is as in the univariate case (see~\eqref{eq:gaus x}).
	\begin{align*}
		~\eqref{term: multigaus x}\le\frac{\ell A}{\theta_{\min}^2}\cdot\norm{\gS_n^{-\half}}_{\op}\frac{K_n^2}{\sqrt{n}}\lesssim D_1 \frac{\log n}{\sqrt{n}}.
        \end{align*}
	Putting this all together and optimizing over $K_n$, similarly to how we did it in the univariate case yields
    \begin{align*}
		\sup_{\mvx} &\abs{\pr\left(\sum_{i=1}^n\mvY_{i,\mvx}\succeq-\sqrt{n}\cdot\mvx_n\right)-\pr\left(\gS_n^{\half}\mvZ\succeq-\mvx_n\right)}\notag\\
		&\lesssim  \ell n^{-1/2}+D_1^{\frac32}\cdot\sqrt{\frac{\log n}{n}}+D_1\frac{\log n}{\sqrt{n}}
         +\frac{\ell^{1/4}}{{\sqrt{n}}}D_1\left(D_2+D_3\ell^{-1}\right),
	\end{align*}
    which holds for $n\ge N(\eps,c,C)$. Multiplying the right-hand side by a universal constant that is large enough yields the bound for all $n$.
\end{proof}

\section{Optimal rate of convergence}\label{sec:optimal}
In the situation of Theorem~\ref{theorem: univ DG med}, suppose that $D_1, D_2,$ and $D_3$ are also bounded by constants, Theorem~\ref{theorem: univ DG med} gives the bound of order $\log n/\sqrt{n}$. In this section, we analyze the i.i.d.~case to motivate the following conjecture about the optimal rate.

\begin{conj}\label{conj:optimal}
	In the situation of Theorem~\ref{theorem: univ DG med}, if $\max_{i\in\{1,2,3\}}\{D_i\}<C$ for some $C>0$, then
	$$\sup_{x\in\dR}\abs{\pr\left(\theta_n\sqrt{n}\cdot M_n\le\gs_nx\right)-\Phi(x)}\approx n^{-\half}.$$
\end{conj}

\begin{lem}
	Let $n=2m+1$ for some integer $m$. Suppose $(X_i)_{i\in[n]}$ are i.i.d.~continuous random variables with CDF $F(x)$ that satisfy Assumption~\ref{as1}.\ref{as: unique med q}--\ref{as: cont q} with $\ga=\half$ and $m_\ga=0$. Assume that $F''(0)$ exists. Then
	$$
	n^{\half}\cdot\sup_{x\in\dR}\abs{\pr\left(\theta\sqrt{2m}\cdot M_n\le x/2\right)-\Phi(x)}\to\frac{|F''(0)|}{4F'(0)^2}\cdot\sup_{x\in\dR}x^{2}\phi(x)=\frac{1}{\sqrt{8\pi e}}\cdot\frac{|F''(0)|}{F'(0)^2}.$$
\end{lem}
Recall that when $X_i$ are uncorrelated we have that $\gs\equiv\half$ and $\theta=F'(0).$
\begin{proof}
	By definition of the sample median and independence of $X_i$'s, we have
	\begin{align}
		\pr\left( 2\theta M_n\le x/\sqrt{2m}\right) & =(2m+1)\binom{2m}{m}\int^{x/(2\theta\sqrt{2m})}_{-\infty}(1-F(t))^mF(t)^m\,dF(t)\notag\\
		& =\frac{2m+1}{2^{2m}}\binom{2m}{m}\int^{x/(2\theta\sqrt{2m})}_{-\infty}(1-(2F(t)-1)^{2})^m\,d F(t)\notag\\
		& =\frac{2m+1}{\sqrt{2m}\cdot 2^{2m+1}}\binom{2m}{m}\int^{x_m}_{-\sqrt{2m}}(1-t^2/2m)^m\,dt,\label{eq:changeu}
	\end{align}
	where $x_m = 2\sqrt{2m}(F(x/(2\theta\sqrt{2m}))-F(0))$. Here, in the last equality, we use the change of variable $t=2\sqrt{2m}\cdot (F(t)-1/2)$. Define
	\[
		a_{m} =\sqrt{2\pi}\cdot\frac{2m+1}{\sqrt{2m}\cdot 2^{2m+1}}\binom{2m}{m}=1+\frac{3}{8m}+O(m^{-2}),
	\]
	where the last equality follows by Stirling's approximation. Thus, we have
	\begin{align*}
		~\eqref{eq:changeu} & =a_{m}\int^{x_m}_{-\sqrt{2m}}\left(g(t^2/2m)\right)^m\phi(t)\,dt
	\end{align*}
	where $g(s):=\left(1-s\right)\exp\left(s\right)$. Notice that for all $\abs{s}\le 1$, we have
	\[
		1-s^{2}\le g(s)\le 1-s^2/4.
	\]

	Thus we get,
	\begin{align*}
		\pr\left(\theta M_n\le x/\sqrt{2m}\right)-\Phi(x)
		 & =\int_{-\sqrt{2m}}^{x_m}\left(g(t^2/2m)^{m}-1\right)\phi(t)\,dt
		+\int_{x}^{x_m}\phi(t)dt+O(1/m)\\
		 &\approx\int_{-\sqrt{2m}}^{x_m}\frac{t^4}{2m}\phi(t)\,dt+\phi(x)(x_m-x)+O(1/m)\\
		 & = O(1/m)+\phi(x)(x_m-x).
	\end{align*}
	Recall that $x_m=2\sqrt{2m}\left(F\left(x/(2\theta\sqrt{2m})\right)-F(0)\right)$.
	Hence, using Taylor expansion, we have,
	\begin{align*}
		x_m-x & =2\sqrt{2m}\left(\frac{x}{2\theta\sqrt{2m}}\cdot F'(0)+\frac{ x^2}{16\theta^2m}F''(0)+R\right)-x=\frac{F''(0)}{4F'(0)^2}\cdot\frac{x^2}{\sqrt{n}}(1+o(1)),
	\end{align*}
	where $R=O(1/n)$. This completes the proof.
\end{proof}

\section*{Acknowledgements}
We would like to thank Sabyasachi Chatterjee for many insightful conversations. We thank the anonymous referees for their careful reading and insightful comments, which resulted in an improved presentation of the article.

\bibliographystyle{alea3}
\bibliography{median}

\end{document}